\documentclass[12pt]{article}
\usepackage{graphicx}
\def\forall{\hbox{for all}~}

\def\L{{\bf L}}

\def\bfw{{\bf w}}
\def\bfu{{\bf u}}
\def\bfq{{\bf q}}

\def\ve{\varepsilon}
\def\n{\noindent}

\def\D{{\cal D}}
\def\N{{\cal N}}

\def\Q{{\cal Q}}

\def\R{I\!\!R}
\def\AC{{\cal AC}}

\def\vp{\varphi}

\def\P{{\cal P}}

\def\vs{\vskip 2em}

\def\bfp{{\bf p}}
\def\v{\vskip 1em}
\def\O{{\cal O}}
\def\C{{\cal C}}
\def\ov{\overline}

\def\Hat{\widehat}

\def\bel{\begin{equation}\label}
\def\eeq{\end{equation}}
\def\sqr#1#2{\vbox{\hrule height .#2pt
\hbox{\vrule width .#2pt height #1pt \kern #1pt
\vrule width .#2pt}\hrule height .#2pt }}
\def\square{\sqr74}
\def\endproof{\hphantom{MM}\hfill\llap{$\square$}\goodbreak}

\begin{document}
 \[\begin{array}{cc}\hbox{\LARGE{\bf On the Controllability of Lagrangian Systems } }
   \\\\ 
 \hbox{\LARGE{\bf by Active Constraints}}\end{array}\]
 
\[\hbox{Alberto Bressan$^\ast$ and Zipeng Wang$^\dagger$}\]
 
 \[\begin{array}{cc}\hbox{($\ast$) Department of Mathematics, Penn State University,  }
   \\
 \hbox{University Park, PA, 16802, USA}\end{array}\]
 \[\hbox{email: Bressan@math.psu.edu}\]
 
  \[\begin{array}{cc}\hbox{($\dagger$) Cambridge Centre for Analysis, University of Cambridge }
   \\
 \hbox{Cambridge, CB3 0WA, UK}\end{array}\]
  \[\hbox{email:  zipeng.wang@cantab.net}\]

\begin{abstract}
We consider a mechanical system which is controlled by means of
moving constraints. Namely, we assume that some of the coordinates can be
directly assigned as functions of time by means of frictionless constraints.
This leads to a system of ODE's whose right hand side depends quadratically
on the time derivative of the control.  In this paper we introduce a
simplified dynamics, described by a differential inclusion.
We prove that every trajectory of the differential inclusion can be uniformly
approximated by a trajectory of the original system, on a sufficiently
large time interval, starting at rest. Under a somewhat stronger assumption,
we show this second trajectory reaches exactly the same terminal point.
\end{abstract}

\section{ Introduction}
\label{sec:1}
\setcounter{equation}{0}
Consider a system whose state is described by $N$ Lagrangian variables
$q^1,\ldots, q^N$.  Let the
 kinetic energy $T=T(q,\dot q)$
be given by a positive definite quadratic form
of the time derivatives $\dot q^i$, namely
\bel{3.1}
T(q,\dot q)~=~\frac{1}{2} \dot q^\dagger G \dot q~=~{1\over 2}
\sum_{i,j=1}^N g_{ij}(q) \,\dot q^i\dot q^j\,.
\eeq
Let the coordinates be split in two groups:
$\{q^1,\ldots, q^n\}$ and $\{ q^{n+1},\ldots, q^{n+m}\}$, with $N=n+m$.
The $(n+m)\times (n+m)$ symmetric
matrix $G$ in (\ref{3.1}) will thus take the corresponding block form
\bel{Gblock}
G=\left(\begin{array}{cc} G_{11}& G_{12}\\ G_{21} &G_{22}\end{array}
\right)=
\left(\begin{array}{cc} (g_{ij})& (g_{i, n+\beta})\\
(g_{n+\alpha, j})  &(g_{n+\alpha,\,n+\beta})
\end{array}
\right).\eeq
We assume that
a controller can prescribe the values
of the last $m$ coordinates
as functions of time, say
\bel{3.3}
q^{n+\alpha}(t)=u_\alpha(t)\qquad\qquad \alpha=1,\ldots,m\,,
\eeq
by implementing $m$ frictionless constraints.
Here {\bf frictionless} means that the forces produced by the constraints
make zero work in connection with any virtual displacement of the
remaining free coordinates $q^1,\ldots,q^n$.
In the absence of external forces, the motion is thus governed by the equations
\bel{lagreq}
\frac{d}{dt}
\frac{\partial T}{\partial \dot q^i} (q,\dot q)-\frac{\partial T}{\partial q^i}
(q,\dot q)~=~ \Phi_i(t)\qquad\qquad i=1,\ldots, n+m\,.
\eeq
Here $\Phi_i$ are the components of the forces generated by the constraints.
The assumption that these constraints are frictionless is expressed by the
identities
\bel{frictionless}
\Phi_i (t)\equiv 0\qquad\qquad i=1,\ldots,n\,.
\eeq
By introducing the conjugate momenta
\bel{conjmom}
p_i~=~p_i(q,\dot q)~\doteq~ {\partial T\over\partial \dot q^i}
~=~ \sum_{i=1}^{n+m} g_{ij}(q)\,\dot q^j\,,
\eeq
it is well known that the evolution of the first $n$ variables
$(q^1,\ldots, q^n)$ and of the corresponding momenta
$(p_1,\ldots, p_n)$ can be described by the system
\bel{Smech}
\left(\begin{array}{c}\dot q\\\,\\\dot p\end{array}\right) =
  \left(\begin{array}{c} Ap\\\, \\-\frac{1}{2}p^\dagger\frac{\partial
A}{\partial
q}p \end{array}\right)
+\left(\begin{array}{c} { K}\\\,\\
-p^\dagger\frac{\partial { K}}{\partial q} \end{array}\right)\dot u\, +\,\,
\dot u^{\dag} \left(\begin{array}{c} {0}\\\,\\
\frac{1}{2} \frac{\partial E}{\partial
q}\end{array}\right)\dot u  \,.
\eeq
Here $A,K,E$ are functions of $q,u$, defined as
\bel{AKE}
A =\left(a^{ij}\right)\doteq (G_{11})^{-1}\,,\qquad\qquad
E=G_{22}- G_{21}AG_{12}\,,
\qquad\qquad
K= -A G_{12}\,.
\eeq
For convenience, in (\ref{Smech}) the vectors $q,p\in\R^n$
are written as column vectors, while the symbol $^\dagger$
denotes transposition.

In general, (1.7) is a system of equations whose right hand side
depends quadratically on the time derivatives of the control
function $u=(u_1,\ldots, u_m)$.
A detailed description of all trajectories of this system is difficult,
because of the interplay between
linear and quadratic terms.  In this paper, to study
(1.7) we introduce a simplified system, described by a differential
inclusion.  For each $q,u$, we define  the convex cone
$$
\Gamma(q,u)\doteq \ov{co} \left\{A(q,u) \Big(w^\dagger
\frac{\partial E(q,u)}{\partial
q}\, w\Big)\,;~~w\in \R^m\right\},
$$
where $\ov{co}$ denotes a closed convex hull.
Intuitively, one can think of $\Gamma(q,u)$ as the set of velocities
which can be instantaneously produced at $(q,u)$, by small vibrations
of the active constraint $u(\cdot)$.
We then consider the differential inclusion
\bel{DI4}
\dot q \in K(q,u)\,\dot u +\Gamma(q,u)\qquad\qquad q(0)=\bar q,\quad
u(0)=\bar u.
\eeq
Trajectories of (\ref{DI4}) will be compared with trajectories of
the original system (\ref{Smech}), with initial data
\bel{indata3}
q(0)=\bar q,\qquad u(0)=\bar u,\qquad p(0)= 0.
\eeq
Our main results show that, for every solution $s\mapsto q^*(s)$
of (\ref{DI4}), say defined for $s\in [0,1]$, there exists
a smooth solution $t\mapsto (q(t), \, p(t))$
of the Cauchy problem (\ref{Smech}), (\ref{indata3}), defined on a suitably
long time interval $[0,T]$, following almost the same path.
Namely, given $\ve>0$, a solution $(q,p)$ of (\ref{Smech}), (\ref{indata3})
can be found such that
\bel{b1}\Big| q(t) - q^*(\psi(t))\Big|  <\ve\,,\qquad\quad
|p(t)|<\ve\qquad\qquad \hbox{for all}~t\in [0,T],
\eeq
for a suitable time rescaling $\psi:[0,T]\mapsto [0,1]$.
Under a somewhat stronger assumption,
the terminal values of the two trajectories can be made equal, namely
$$q(T) = q^*(1).$$
{\bf Remark 1.} Since the components $p_i$ bear a linear relation
to the velocities $\dot q^j$, the system (\ref{Smech}) describes a
``second order" dynamics, which could be equivalently written in
terms of the second derivatives $\ddot q^j$.  On the other hand, the
reduced system (\ref{DI4}) contains no inertial term, and is
essentially of first order. The inequalities (\ref{b1}) show that,
keeping $p(t)\approx 0$, the two dynamics can be related.  We remark that the
present results are entirely different in nature from those in \cite{B08, B-R3,
B-R4},
where the impulsive control system is approximated by a differential inclusion
living in the $2n$-dimensional space described by the $(q,p)$-variables.

The
paper is organized as follows. Section 2 contains precise statements
of the main results. The proofs are then worked out in Sections
3--5.  Section 6 contains two examples.  The first one shows
the necessity of a technical assumption. The second one provides a
simple application to the control of a bead sliding without friction
along a rotating bar. The last section is the derivation of  evolution equations in (1.7).

For the theory of multifunctions and
differential inclusions we refer to \cite{AC} or \cite{Smi}.
Earlier results on impulsive control systems were provided in
\cite{B-R1, B-R3, AB1, AB2}. A general introduction to the theory of
control can be found in \cite{BP, BL, NS} and in \cite{So}.  
We remark that the idea of averaging,
used in the proof of our main theorem, is widespread
in the analysis of mechanical systems with oscillatory behavior. 
Several results in this direction can be found 
in \cite{AKN,BoMi}.

\section{Statement of Main Results}
\setcounter{equation}{0}
Motivated by the model (\ref{Smech}), from now on
we consider a system of the form
\bel{DE1}
\left(\begin{array}{c}\dot q\\\,\\\dot p\end{array}\right) =
  \left(\begin{array}{c} A p\\\, \\-\frac{1}{2}p^\dagger B p\end{array}\right)
+\left(\begin{array}{c} { K}\\\,\\
-p^\dagger C  \end{array}\right)\dot u\, +\,\,
\dot u^{\dag} \left(\begin{array}{c} {0}\\\,\\
 D \end{array}\right)\dot u  \,.
\eeq
Given an initial data
\bel{indata}
q(0)= \bar q,  \qquad\quad p(0)=\bar p,\quad\qquad u(0)=\bar u,
\eeq
we shall study the set of trajectories
of (\ref{DE1}).

The difficulty in analyzing (\ref{DE1}) stems from the fact that the right
hand side
contains both linear and
quadratic terms w.r.t.~the time derivative $\dot u$.
A simplification can be achieved by considering separately the contributions of
these terms.
If $D\equiv 0$, we have
the reduced system
\bel{DE2}
\left(\begin{array}{c}\dot q\\\,\\\dot p\end{array}\right) =
  \left(\begin{array}{c} A p\\\, \\-\frac{1}{2}p^\dagger B p\end{array}\right)
+\left(\begin{array}{c} { K}\\\,\\
-p^\dagger C  \end{array}\right)\dot u  \,.
\eeq
Notice that, if $\bar p=0$, then $p(t)\equiv 0$ for every time $t$.
In this case, the trajectory of the system (\ref{DE2}) is entirely determined
by
solving the reduced equation
\bel{DE3}
\dot q= K(q,u) \dot u\,,\qquad\qquad q(0)=\bar q\,.
\eeq
We claim that, even in the case $D\not= 0$,
given a sufficiently long time interval, every trajectory of
(\ref{DE3}) can be uniformly approximated by a trajectory
of the original system (\ref{DE1}).   More generally,
if the initial speed is sufficiently small, then  one can track
every solution to the differential inclusion
\bel{DI1}
\dot q\in K(q,u) \dot u+\Gamma(q,u)\,,\qquad\qquad q(0)=\bar q\,.
\eeq
Here $\Gamma$ is the cone defined by
\bel{Gammadef}
\Gamma(q,u)\doteq \ov{co} \Big\{A(q,u)\, (w^\dagger D(q,u)\, w)\,;
~~w\in \R^m\Big\},
\eeq
where $\ov{co}$ denotes the closed convex hull of a set.

\n{\bf Definition 1.} Given an absolutely continuous control
function $t\mapsto u(t)$,  defined for $t\in [0,T]$,
by a Carath\'eodory solution of the differential inclusion  (\ref{DI1})
we mean an absolutely continuous map $t\mapsto q(t)$
such that
\bel{DI2}
\dot q(t)-K(q(t), u(t)) \dot u(t)\in \Gamma(q(t), u(t))
\qquad\qquad\hbox{for a.e.}
\quad t\in [0,T].
\eeq
Our main result is concerned with approximation of trajectories of
(\ref{DI1}) with solutions of the full system (\ref{DE1}).
Our basic hypotheses are as follows.

\n{\bf (H)} The matrices $A, B, K, C$ in (\ref{DE1}) are locally
Lipschitz continuous functions of the variables $q,u$, and the same
is true of $D$ and of the partial derivatives $K_q, K_u$. Moreover,
the cone $\Gamma$ in (\ref{Gammadef}) depends continuously on
$(q,u)$; namely, the compact, convex valued multifunction
\bel{Gamma1} (q,u)~\mapsto~ \Gamma_1(q,u)~\doteq~ \Big\{ p\in
\Gamma(q,u)\,;~~|p|\leq 1\Big\} \eeq
 is continuous w.r.t.~the
Hausdorff distance. 

{\bf Theorem 1.} {\it  Let the assumptions
(H) hold, and let $s\mapsto q^*(s)$ be any Carath\'eodory solution
of differential inclusion (\ref{DI1}) defined for $s\in [0,1]$,
corresponding to an absolutely  continuous control $u^*(\cdot)$.

Then, for every $\ve>0$, there exists $\delta>0$, an interval $[0,T]$ and a
smooth control $u(\cdot)$
defined on $[0,T]$ such that the following holds.
If $|\bar p|<\delta$, then the corresponding solution of (\ref{DE1})
with initial data
(\ref{indata})
satisfies
\bel{approx2}\sup_{t\in [0,T]} |p(t)| <\ve,\qquad\quad
\sup_{t\in [0,T]} \Big|q(t)-q^*(\psi(t))\Big| <\ve,\quad\qquad \sup_{t\in
[0,T]}
\Big|u(t)-u^*(\psi(t))\Big| <\ve,
\eeq
for some increasing diffeomorphism $\psi: [0,T]\mapsto [0,1]$.}

{\bf Remark 2.}   Assume that, more generally,
the control $u^*$ and the trajectory $q^*$
are defined on an interval $[0,T^*]$.
Since $\dot u$ enters linearly in the equation (\ref{DI1}),
 the rescaled function $ \hat q(s)=q^*(T^*s)$
provides another solution of (\ref{DI1}), corresponding to the control
$\hat u(s)\doteq u^*(T^*s)$.
By a linear rescaling of time , it is thus not restrictive to assume that
$q^*, u^*$ are defined for $s\in [0,1]$.

Next, we consider the problem of exactly reaching a state
$(\bfq,\bfu)$ at some (possibly large) time $T$, with small terminal
speed. As a preliminary, we introduce a notion of normal
reachability. As in \cite{J}, this means that there exists a family
of trajectories whose terminal points nicely cover a whole
neighborhood of the target point $(\bfq,\bfu)$. More precisely: \v
{\bf Definition 2.}  Given the differential inclusion (\ref{DI1}),
the state $(\bfq,\bfu)\in R^{n+m}$ is {\em normally reachable} from
the initial state $(\bar q,\bar u)$ if there exists a parameterized
family of trajectories
$$s\mapsto (q^\lambda(s), \,u^\lambda(s)), \qquad\qquad
\lambda\in \Lambda\subset\R^{n+m}, ~~s\in [0,1]$$
with the following properties.

(i)  The parameter $\lambda$ ranges in a neighborhood
$\Lambda$ of the origin in $\R^{n+m}$.  The map $\lambda\mapsto
(q^\lambda(\cdot),\,u^\lambda(\cdot))$ is continuous from $\Lambda$ into
$W^{1,1}\left([0,1]\,;~\R^{n+m}\right)$.

(ii) For every $\lambda\in\Lambda$ we have $(q^\lambda(0),\,u^\lambda(0))
=(\bar q,\bar u)$.   Moreover, when $\lambda=0\in\R^{n+m}$ we have
$(q^0(1),\,u^0(1))= (\bfq,\bfu)$ and
the $(n+m)\times (n+m)$ Jacobian matrix
$$\left({\partial(q^\lambda(1),\,u^\lambda(1))\over
\partial\lambda}\right)$$
has full rank, i.e.~it is invertible.

{\bf Theorem 2.} {\it Let (H) hold, and assume that the state
$(\bfq,\bfu)\in R^{n+m}$ is {\em normally reachable} from the
initial state $(\bar q,\bar u)$, for the differential inclusion
(\ref{DI1}).   Then, for any $\ve>0$, there exists $\delta>0$ such
that the following holds. If $|\bar p|<\delta$, there exists a time
$T$ and a control function $u$ defined on $[0,T]$  such that the
corresponding solution of (\ref{DE1}) satisfies (\ref{approx2})
together with \bel{exact} (q(T),\,u(T)) =(\bfq,\bfu)\,.\eeq }

The proof of Theorem 2 relies on a topological argument.
The key ingredient is the following continuous approximation lemma.
By $\AC([0,T])$ we denote here
the space of absolutely continuous functions on $[0,T]$, with norm
\[\|f\|_{\AC}~\doteq ~\int_0^T|\dot f(t)|\,dt+\sup_{t\in [0,T]}~|f(t)| \,.\]
{\bf Lemma 1.} {\it Let (H) hold.
Consider a family of solutions
$(q^\lambda, u^\lambda)$ of the differential inclusion (\ref{DI1}),
assuming that the map $\lambda\mapsto
(q^\lambda(\cdot),\,u^\lambda(\cdot))$ is continuous from
a compact set $\Lambda\subset\R^d$ into
$\AC\left([0,1]\,;~\R^{n+m}\right)$.
Then, given any $\ve>0$, there exists a map
$(\lambda, s)\mapsto (\tilde u^\lambda(s), \,\tilde  w^\lambda(s))$
from $\Lambda \times [0,1]$ into $\R^m\times \R^m$,
which is continuous w.r.t.~$\lambda$
and $\C^\infty$ in the variable $s$,
such that the following holds.
Calling $\tilde q^\lambda(\cdot)$ the solution to
\bel{qapp}
{d\over ds} q(s)~= ~K(q(s),\tilde
u^\lambda(s))\cdot {d\over ds} \tilde u^\lambda(s)+A(q(s),
\tilde u^\lambda(s))\Big(\tilde w^\lambda(s)^\dagger D(q(s), \tilde u^\lambda
(s))
\tilde w^\lambda(s)\Big),\qquad\quad q(0)=\bar q,
\eeq
for every $\lambda\in\Lambda$ one has
\bel{app5}
\sup_{s\in [0,1]} \,\left|
\tilde q^\lambda(s)-q^\lambda(s)\right| <\ve,\quad\qquad
\sup_{s\in [0,1]}\,
|\tilde u^\lambda(s)-u^\lambda(s)| <\ve.
\eeq
}

{\bf Remark 3.} The assumption (H) requires
that the maps $A, B, K, C, D$
be locally Lipschitz continuous.  We observe that, toward the proof of
Lemma 1, it is not restrictive to assume that
all these maps have compact support,
and  are therefore globally Lipschitz continuous.  Indeed,
the set
\bel{Omega0}\Omega_0\doteq \Big\{ (q^\lambda(s),\, u^\lambda(s))\,;~~
s\in [0,1]\,,~~\lambda\in\Lambda\Big\}\,,
\eeq
is compact, and the same is true for its closed neighborhood
\bel{Omegarho}\Omega_\rho\doteq\Big\{ (q,u)\in\R^{n+m}\,;~~
|q-q^\lambda(s)|\leq \rho\,,~~
|u-u^\lambda(s)|\leq \rho\qquad\hbox{for some}~~
s\in [0,1]\,,~~\lambda\in\Lambda\Big\},
\eeq
for any $\rho>0$.
Let $\varphi:\R^{n+m}\mapsto \R_+$ be a smooth cutoff function
such that
$$\varphi(q,u)~=\left\{ \begin{array}{rl} 1&\quad\hbox{if}~~(q,u)\in
\Omega_1\,,
\cr
0&\quad\hbox{if}~~(q,u)\notin \Omega_2\,.
\end{array}
\right.
$$
The functions $\Hat A\doteq \vp\cdot A, \ldots,\Hat D\doteq \vp\cdot
D$ have compact support and are thus globally Lipschitz continuous.
If the conclusion of Lemma 1 holds for $\Hat A,\Hat B, \Hat K,\Hat
C,\Hat D$, then it also holds for the original functions $A, B, K,
C, D$. Indeed, when $0<\ve<1$, the inequalities (\ref{app5}) imply
that $(q^\lambda(s), p^\lambda(s))\in \Omega_1$. Restricted to
$\Omega_1$, one has the identities $A=\Hat A,\ldots, D=\Hat D$. This
same remark applies to the proofs of Theorems 1 and 2.

\section{Proof of the approximation lemma}
\label{sec:2}
\setcounter{equation}{0}
We first prove two auxiliary results. Recall that the convex sets
$\Gamma_1(q,u)$ were
defined at (\ref{Gamma1}).
For notational convenience, we introduce the set of coefficients
of convex combinations
$$\Delta_\nu\doteq \Big\{ \theta=(\theta_1,\ldots,\theta_\nu)\,;
~~~\theta_i\in [0,1]\,,~~~
\sum_{i=1}^\nu\theta_i=1\Big\}.$$
{\bf Lemma 2.}  {\it
Given $\ve'>0$ and a compact set $\Omega\subset\R^{n+m}$, there
exist finitely many vectors $\bfw_1,\ldots, \bfw_\nu$
such that the following holds.
Given any $(q,u)\in \Omega$ and any
$p\in \Gamma_1(q,u)$,
there exist coefficients $(\theta_1,\ldots,\theta_\nu)\in\Delta_\nu$
 such that
\bel{appr6}
\left|p-A(q,u)\cdot\sum_{i=1}^\nu
\theta_i \bfw_i^\dagger D(q,u)\bfw_i\right|\leq\ve'\,.
\eeq
}

{\bf Proof.} Consider the  domain
\bel{domTheta}
\D\doteq \Big\{(q,u,p)\,;~~(q,u)\in \Omega\,,~~p\in \Gamma_1(q,u)\Big\}\,.
\eeq
Notice that $\D$ is compact, because of the assumption (H).
For each $(\bar q,\bar u,\bar p)\in \D$, choose finitely many
vectors $\bfw_i=\bfw_i^{(\bar q,\bar u,\bar p)}$ and coefficients
$\theta_i=\theta_i^{(\bar q,\bar u,\bar p)}$, $i=1,\ldots,M$, such that
$$
\sum_{i=1}^M \theta_i =1\,,\qquad\qquad \left|\bar p-A(\bar q,\bar u)
\cdot\sum_{i=1}^M
\theta_i \bfw_i^\dagger D(\bar q,\bar u)\bfw_i\right|\leq
{\ve'\over 2}\,.
$$
By continuity, we still have
$$
\left|p-A(q,u)\cdot\sum_{i=1}^M
\theta_i \bfw_i^\dagger D(q,u)\bfw_i\right|\leq\ve'
$$
for all $(q,u,p)$ in a neighborhood $V^{(\bar q,\bar u,\bar p)} $ of the point
$(\bar q,\bar u,\bar p)$.    Covering the compact domain $\D$ with finitely
many
neighborhoods $V_\ell=V^{(q_\ell, u_\ell, p_\ell)}$, $\ell=1,\ldots,\kappa$,
and choosing
$$\{\bfw_1,\ldots, \bfw_\nu\}~\doteq ~\Big\{
\bfw_1^{(q_1, u_1, p_1)},\ldots, \bfw_{M(1)}^{( q_1, u_1, p_1)}, ~~~\ldots~~~,
\bfw_1^{( q_\kappa, u_\kappa, p_\kappa)},\ldots,
\bfw_{M(\kappa)}^{( q_\kappa, u_\kappa, p_\kappa)}
\Big\},$$
we achieve the conclusion of the lemma.
\endproof

The next lemma states that, if we relax the inequality in
(\ref{appr6}), the coefficients $\theta_i$ can be chosen
depending continuously on
$q,u,p$.
\v
{\bf Lemma 3.}
{\it
Given a compact set $\Omega\subset\R^{n+m}$, define the compact domain
$\D$ as in (\ref{domTheta}).
Then, for any $\ve'>0$,
there exists a continuous mapping
$\Theta= (\Theta_1,\ldots,\Theta_\nu):\D\mapsto \Delta_\nu$,
such that
\bel{Thapp}
\left|p-A(q,u)\cdot\sum_{i=1}^\nu
 \Theta_i(q,u,p) \bfw_i^\dagger D(q,u)\bfw_i\right|\leq 2\ve'\,.
\eeq
for all $(q,u,p)\in\D$.
}
\v
{\bf Proof.}
By continuity and compactness, there exists $\delta>0$
such that the following holds.
If
\bel{a5}
|q-\tilde q|<\delta\,,\qquad |u-\tilde u|<\delta\,,\qquad |p-\tilde p|<\delta,
\qquad |\theta_i-\tilde\theta_i|<\delta\qquad\hbox{for} ~i=1,\ldots,\nu\,,
\eeq
and if
\bel{a6}
\left|\tilde p-A(\tilde q,\tilde u)\cdot\sum_{i=1}^\nu
\tilde \theta_i \bfw_i^\dagger D(\tilde q,\tilde u)\bfw_i\right|\leq \ve'\,,
\eeq
then
\bel{a7}\left| p-A(q,u)\cdot\sum_{i=1}^\nu
\theta_i \bfw_i^\dagger D( q, u)\bfw_i\right|\leq 2\ve'\,.
\eeq
Next, consider the set-valued function
$$\Hat\Theta(q,u,p)\doteq
\left\{ \theta=(\theta_1,\ldots,\theta_\nu)\in\Delta_\nu\,,
\qquad
\left|p-A(q,u)\cdot\sum_{i=1}^\nu \theta_i\,
\bfw_i^\dagger D(q,u)\bfw_i\right|\leq\ve'\,\right\}.
$$
Observe that the multifunction
$\Hat\Theta:\D\mapsto\Delta_\nu$ has closed graph, and non-empty,
compact, convex values.   By a selection theorem in
\cite{AC}, for every $\delta>0$, this multifunction admits a continuous,
$\delta$-approximate selection $\Theta:\D\mapsto \Delta_\nu$,
in the sense of graph. Calling $\N(S,\delta)$ the $\delta$-neighborhood
around a set $S$, this means that
$$\hbox{graph}\, \Theta ~\subset~\N\left( \hbox{graph}\,
\Hat\Theta\,,\,\delta\right)\,.$$
If $\delta>0$ was chosen sufficiently small,
so that (\ref{a5})-(\ref{a6}) imply (\ref{a7}),
then the continuous function $\Theta$
satisfies the conclusion of the lemma.
\endproof

{\bf Proof of Lemma 1.}
According to Remark 3, we can assume that all functions
$A,B,K,C,D$ have compact support, hence they are all
globally Lipschitz continuous and uniformly bounded.
The proof of the continuous approximation lemma will be given in several steps.

{\bf 1.}
By assumption, for every $\lambda\in\Lambda$ we have
\bel{qld}
\dot q^\lambda(s)=K\left( q^\lambda(s),u^\lambda(s)\right)
\dot u^\lambda(s) + \gamma^\lambda(s)\,,\eeq
where $s\mapsto \gamma^\lambda(s)\in
\Gamma(q^\lambda(s),
\, u^\lambda(s))$ is some measurable map, depending continuously
on $\lambda$ in the $\L^1$ norm.

We claim that it is not restrictive to assume that
the functions $\dot q^\lambda(\cdot)$, $\dot u^\lambda(\cdot)$, and
$\gamma^\lambda(\cdot)$ are uniformly bounded.
Indeed, fix an integer $N$ and define the times $s_i\doteq i/N$.
For each $\lambda$, consider the time rescaling
\bel{tlambda}
t^\lambda(s) \doteq s_{i-1}+{\displaystyle
\int_{s_{i-1}}^s \left(|\dot u^\lambda(t)|+|\gamma^\lambda(t)|+N^{-1}
\right)\,dt\over\displaystyle N\cdot
\int_{s_{i-1}}^{s_i} \left(|\dot u^\lambda(t)|+|\gamma^\lambda(t)|+N^{-1}
\right)\,dt}
\qquad\qquad \hbox{if}\quad s\in [s_{i-1},\, s_i].
\eeq
Observe that the map $s\mapsto t^\lambda(s)$ is strictly increasing,
 satisfies
$$t^\lambda(s_i)= s_i\qquad\qquad \hbox{for all} ~i=0,1,\ldots,N,$$
and has a Lipschitz continuous inverse which we denote by
$t\mapsto s^\lambda(t)$.
We now define
$$q_N^\lambda(t) \doteq q^\lambda(s^\lambda(t)),\qquad
u_N^\lambda(t) \doteq u^\lambda(s^\lambda(t)),
\qquad
\gamma_N^\lambda(t) \doteq \gamma^\lambda(s^\lambda(t))\cdot
\left({d
\over dt} s^\lambda(t)\right).$$
By (\ref{tlambda}), the above definitions yield
\bel{qnldot}
{d\over dt} q^\lambda_N(t) = K\left(q^\lambda_N(t), u^\lambda_N(t)\right)
\cdot {d\over dt} u^\lambda_N(t)
+\gamma^\lambda_N(t).
\eeq
Moreover, for a.e.~$t\in [s_{i-1}, s_i]$, (\ref{tlambda}) implies
$$\left|{d\over dt} u^\lambda_N(t)\right| + |\gamma_N^\lambda(t)|
\leq
N\cdot
\int_{s_{i-1}}^{s_i} \left(|\dot u^\lambda(t)|+|\gamma^\lambda(t)|+N^{-1}
\right)\,dt\,,$$
showing that $\dot u^\lambda_N$ and $\gamma^\lambda_N$ remain
uniformly bounded. The continuity w.r.t.~the parameter $\lambda$
implies that these bounds are uniform as $\lambda$ ranges in the compact set
$\Lambda$. Moreover, the maps $\lambda\mapsto q^\lambda_N(\cdot)$
 and $\lambda\mapsto u^\lambda_N(\cdot)$ are continuous
 from $\Lambda$ into $\AC([0,1])$.

 Finally, for any given $\ve>0$, by choosing the integer $N$
sufficiently large we can achieve the inequalities
\bel{appl5}
\sup_{s\in [0,1]} \,\left|
q^\lambda_N(s)-q^\lambda(s)\right| <\ve,\quad\qquad
\sup_{s\in [0,1]}\,
\left|u^\lambda_N(s)-u^\lambda(s)\right| <\ve.
\eeq
Since $q_N^\lambda$ satisfies
(\ref{qnldot}) and $K$ is bounded,
we conclude that the derivative
$\dot q^\lambda_N$ is uniformly bounded as well.
This completes the proof of our claim.
\v
{\bf 2.}
From now on, we can thus assume that
\bel{qugbound}
|\dot q^\lambda(s)|+|\dot u^\lambda(s)|+|\gamma^\lambda(s)|~\leq~ \ov M
\qquad\qquad \hbox{for a.e.}~s\in [0,1],
\eeq
for some constant $\ov M$ and  every $\lambda\in \Lambda$.

Consider the compact set
$\Omega\doteq \Omega_1$
defined as in (\ref{Omegarho}),
and the corresponding domain $\D$ as in (\ref{domTheta}).

For a given $\ve'>0$, whose precise value will
be determined later, we can
choose vectors $\bfw_1,\ldots,\bfw_\nu$
according to Lemma 2.  Let $\Theta:\D\mapsto \Delta_\nu$
be the continuous map constructed in Lemma 3, and
define the measurable coefficients
\bel{thdef}
\theta_i^\lambda(s)\doteq \Theta\left(q^\lambda(s),\, u^\lambda(s),\,
{\gamma^\lambda(s)\over |\gamma^\lambda(s)|}\right).
\eeq
By (\ref{Thapp}) we have
\bel{appr7}
\left|\gamma^\lambda(s)-A\left(q^\lambda(s),\,u^\lambda(s)\right)
\cdot\sum_{i=1}^\nu |\gamma^\lambda(s)|\,
\theta_i^\lambda(s)\bfw_i^\dagger D
\left(q^\lambda(s),\,u^\lambda(s)\right)\bfw_i\right|
~\leq~ 2\ve'\,|\gamma^\lambda(s)|~\leq~  2\ve'\ov M
\eeq
for a.e.~$s\in [0,1]$.

{\bf 3.} Next, we divide the interval $[0,1]$ into $k\,\nu$  equal
subintervals,
choosing $k$ very large.  For notational convenience we set
$$\tau_j\doteq {j\over k}\,,\qquad\qquad \tau_{j,\ell}\doteq {j\over k} +
{\ell\over k\,\nu}\,.$$
Here $j=0,\ldots,k$, while $\ell=0,\ldots,\nu$.
For each $\lambda\in\Lambda$, we now define a
continuous,  piecewise affine control function
$s\mapsto\tilde u^\lambda(s)$ by setting
\bel{cudef}
\tilde u^\lambda(\tau_j)~\doteq~ u^\lambda(\tau_j) \qquad \qquad
j=0,\ldots, k,
\eeq
and extending $\tilde u^\lambda$ to an affine map on each interval
$\left[\tau_{j-1}\,,~\tau_{j} \right]$.
Since by (\ref{qugbound})
the functions $u^\lambda(\cdot)$ are uniformly
Lipschitz continuous, by choosing $k$  large enough we can  achieve
the bounds
\bel{Qd4}
\left|\tilde u^\lambda(s)-u^\lambda(s)\right|<\ve\qquad\qquad
\forall s\in [0,1],~~\lambda\in\Lambda.
\eeq
Moreover, we define
\bel{cwdef}
\tilde w^\lambda(s)~\doteq~ \left(k\,\nu\cdot \int_{\tau_j}^{\tau_{j+1}}
|\gamma^\lambda(s)|\,
\theta^\lambda_\ell(s)\, ds\right)^{1/2}\bfw_\ell\qquad\qquad
\hbox{for}\qquad s\in \,]\tau_{j,\ell-1}, ~\tau_{j,\ell}]\,.
\eeq
Here $j=0,\ldots,k-1$, while $\ell=1,\ldots,\nu$.
Call $\tilde q^\lambda(\cdot)$
the corresponding solution
of (\ref{qapp}). In the next step we will prove that,
by choosing first $\ve'>0$ sufficiently small and
then the integer $k$ large enough, the inequalities in
(\ref{app5}) are satisfied.

{\bf 4.}
To compare the two functions $q^\lambda(\cdot)$ and
$\tilde q^\lambda(\cdot)$, we introduce a third function $Q^\lambda(\cdot)$,
defined as the solution to the
Cauchy problem
\bel{cp3}\begin{array}{rl}
&\displaystyle\dot Q^\lambda(s)=K\left( Q^\lambda(s),u^\lambda(s)\right)
\dot u^\lambda(s) +A\left( Q^\lambda(s),u^\lambda(s)\right)\cdot
\sum_{i=1}^\nu
|\gamma^\lambda(s)|\,\theta_i^\lambda(s)\bfw_i^\dagger
D\left(Q^\lambda(s),\,u^\lambda(s)\right)\bfw_i\,,\cr
& Q^\lambda
(0)=\bar q\,.
\end{array}
\eeq
To estimate the difference $|q^\lambda-Q^\lambda|$,
consider the Picard map $y\mapsto \P y$~ (depending on $\lambda\in\Lambda$),
defined as
\bel{picard}
(\P y)(t)\doteq \bar q+\int_0^t
\left\{K(y(s), u^\lambda(s))\,\dot u^\lambda(s)
+\gamma^\lambda(s)\right\}
\, ds\,.
\eeq
By Remark 3 and by step {\bf 1} in this proof,
we can assume that $K$ is globally Lipschitz continuous and that the
functions $\dot u^\lambda$ are uniformly
bounded.  Therefore there exists a constant $L$, independent of
$\lambda\in\Lambda$, such that each Picard
map $\P$ is a strict contraction w.r.t.~the weighted
norm
$$\|y\|_*\doteq \sup_{s\in [0,1]} \, e^{-Ls} |y(s)|\,.$$
More precisely, for every continuous functions $y,\tilde y$,
\bel{scontraction}\|\P y-\P \tilde y\|_*\leq
{1\over 2}\|y-\tilde y\|_*\,.
\eeq
In turn (see for example the Appendix in \cite{BP}),
since $ q^\lambda(\cdot)$ is the fixed point of $\P$,
for every $y(\cdot)$ this implies the estimate
\bel{yqest}
\|y-q^\lambda\|_*\leq 2\|y-\P y\|_*\,.
\eeq
We now have
$$(\P Q^\lambda)(t)-Q^\lambda(t)
=\int_0^t\left\{\gamma^\lambda(s)-A\left(q^\lambda(s),\,u^\lambda(s)\right)
\cdot\sum_{i=1}^\nu |\gamma^\lambda(s)|\,
\theta_i^\lambda(s)\bfw_i^\dagger D
\left(q^\lambda(s),\,u^\lambda(s)\right)\bfw_i\right\}\, ds\,.
$$
By (\ref{appr7}), this yields
\bel{Qd1}\sup_{t\in [0,1]} \left| q^\lambda(t)-Q^\lambda(t)\right|
~\leq~e^L\|q^\lambda-Q^\lambda\|_*~\leq~2e^L\|q^\lambda-\P
q^\lambda\|_*~\leq~4 e^L \ve' \ov M\,.\eeq
Notice that the constant $L$ depends only on the Lipschitz norm of $K$
 and on the upper bound  on $|\dot u^\lambda|$
at (\ref{qugbound}).  Therefore, we can assume that $\ve'>0$
in (\ref{qugbound}) was chosen
so that
\bel{Qd2}\sup_{t\in [0,1]}\left| q^\lambda(t)-Q^\lambda(t)\right|~\leq~
4 e^L \ve' \ov M~<~ {\ve\over 2}\,.
\eeq
Next, to estimate the difference
$|\tilde q^\lambda-Q^\lambda|$, we consider a
second Picard map $y\mapsto \P y$, with
\bel{pic2}
(\P y)(t)\doteq \bar q+\int_0^t
\left\{K(y(s), u^\lambda(s))\dot u^\lambda(s)
+A(y(s), u^\lambda(s))\cdot
\sum_{i=1}^\nu
|\gamma^\lambda(s)|\,\theta_i^\lambda(s)\bfw_i^\dagger
D(y(s),\,u^\lambda(s))\bfw_i\right\}
\, ds
\eeq
By the boundedness of
$\dot u^\lambda$, $\gamma^\lambda$, and by the Lipschitz
continuity of $K,A,D$, this map will be a strict contraction
and satisfy (\ref{scontraction})
w.r.t.~some weighted norm of the form
$$\|y\|_*\doteq \sup_{s\in [0,1]} \, e^{-L's} |y(s)|\,.$$
Notice that in this case the constant $L'$ may depend also on
$\max\{ |\bfw_1|,\ldots,|\bfw_\nu|\}$, and hence on the earlier
choice of $\ve'$.

In addition to (\ref{cudef}), for every $j$ and every choice of the
constants $y_j,u_j$, the definition (\ref{cwdef}) yields
$$\int_{\tau_{j}}^{\tau_{j+1}} \left[\tilde w^\lambda(s)
\right]^\dagger D(y_j, u_j)\, \tilde w^\lambda(s)\,ds
~=~
\int_{\tau_j}^{\tau_{j+1}}\sum_{\ell=1}^\nu
|\gamma^\lambda(s)|\,
\theta^\lambda_\ell(s)\, \bfw_\ell^\dagger\,
D(y_j, u_j)\,\bfw_\ell\,ds.
$$
Therefore, by the uniform Lipschitz continuity of the maps $K,A,D$ and
$\tilde q^\lambda, u^\lambda$, it follows
the estimate
$$
\Big|(\P\tilde q^\lambda)(\tau_j)-\tilde q^\lambda(\tau_j)\Big|
~\leq~
C_1\cdot \sup_{i=1,\ldots,j} ~~\sup_{t,t'\in [\tau_{i-1}, \tau_i]}\left(
\left|\tilde q^\lambda(t)-\tilde q^\lambda(t')\right|+
\left|u^\lambda(t)-u^\lambda(t')\right|\right) ~\leq~ {C_2\over k}\,,
$$
for suitable constants $C_1,C_2$, depending on $\ov M$ and
$\max\{ |\bfw_1|,\ldots,|\bfw_\nu|\}$
but not on $\lambda, k$.
More generally, for $t\in [\tau_{j}, \tau_{j+1}]$ we have
$$
\Big|(\P\tilde q^\lambda)(t)-\tilde q^\lambda(t)\Big|
~\leq~ \Big|(\P\tilde q^\lambda)(t)-(\P\tilde q^\lambda)(\tau_j)\Big|
+\Big|(\P\tilde q^\lambda)(\tau_j)-\tilde q^\lambda(\tau_j)\Big|
+\Big|\tilde q^\lambda(\tau_j)-\tilde q^\lambda(t)\Big|~\leq~
{C_3\over k}\,,
$$
for a suitable constant $C_3$.
Observing that $Q^\lambda$ is the fixed point of the Picard map $\P$
in
(\ref{pic2}),
we can thus choose $k$ large enough so that
\bel{Qd3}
\sup_{t\in [0,1]} \left| \tilde q^\lambda(t)-Q^\lambda(t)\right|
~\leq~e^{L'}\|\tilde q^\lambda-Q^\lambda\|_*~
\leq~2e^{L'}\|\tilde q^\lambda-\P
\tilde q^\lambda\|_*~\leq~2 e^{L'} \cdot {C_3\over k}~<~ {\ve\over 2}\,.
\eeq

{\bf 5.} At this stage we have constructed functions
$\tilde u^\lambda,\tilde w^\lambda$ which satisfy (\ref{app5}).
However, the maps
$$(\lambda, s)\mapsto {d\over ds} \tilde u^\lambda(s),\qquad\qquad
(\lambda,s)\mapsto
\tilde w^\lambda(s)$$
are continuous as functions of $\lambda$, but piecewise constant
with jumps at the points
$\tau_{j,\ell}$ as functions of the time variable $s\in [0,1]$.
To complete the proof, we need to achieve smoothness w.r.t.~the variable $s$.
This is obtained by a standard mollification procedure.

We first extend each the functions  $\tilde u^\lambda$ by setting
$\tilde u^\lambda(s) = \tilde u^\lambda(0)$ if $s<0$,
$\tilde u^\lambda(s) = \tilde u^\lambda(1)$ if $s>1$, and similarly
for $\tilde w^\lambda$.
Then we perform a  mollification in the $s$-variable:
$$ U^\lambda(s)\doteq \int u^\lambda(s-\sigma)\phi_\rho(\sigma)\,d\sigma\,,
\qquad
 W^\lambda(s)\doteq \int w^\lambda(s-\sigma)\phi_\rho(\sigma)\,d\sigma\,.
$$
Here  $\phi_\rho$ is a standard mollification kernel, so that
$\phi_\rho(\sigma)\doteq \rho^{-1}\phi(\rho^{-1} \sigma)$
for some smooth
function with compact support $\phi\in \C^\infty_c\,$, with
$\phi\geq 0$ and $\int\phi(\sigma)\,d\sigma=1$.

By choosing $\rho>0$ sufficiently small, it is clear that the functions
$ U^\lambda$ and
$ W^\lambda$, in place of $\tilde u^\lambda$ and
$\tilde q^\lambda$,  satisfy all conclusions of  Lemma 1.
\endproof

{\bf Remark 4.} Since the solution of (\ref{qapp}) depends continuously
on $\tilde w^\lambda$, we can slightly perturb these functions
in $\L^1$ and still achieve the pointwise inequalities (\ref{app5}).
In particular, on the smooth functions $\tilde w^\lambda$ we
can impose   the additional requirement that
\bel{wvanish}
\tilde w^\lambda(s) =0\qquad\qquad \forall \lambda\in\Lambda\,,\quad
s\in [0,\ve_0]\,,\eeq
for some $\ve_0>0$ sufficiently small.

\section{Proof of Theorem 1}
\label{sec:3}
\setcounter{equation}{0}
Using Lemma 1 in the special case where the parameter set
$\Lambda$ is a singleton,
we can assume that $q^*(\cdot)$ and $u^*(\cdot)$ are smooth,
and that there exists a smooth function $w^*(\cdot)$
such that
\bel{qapp2}
\dot q^*(s)~= ~K(q^*(s),
u^*(s)) \dot u^*(s)+A(q^*(s),
u^*(s))\Big( w^*(s)^\dagger D(q^*(s), u^*(s))
w^*(s)\Big).
\eeq
Moreover, by Remark 4, for some $\ve_0>0$ sufficiently small we can assume that
\bel{wvan}
w^*(s)=0\qquad\qquad \forall s\in [0,\ve_0].\eeq

Define the nonlinear time rescaling $\psi:[0,T]\mapsto [0,1]$,
\bel{psidef}s~=~\psi(t)~\doteq~ 1-{\ln(1+T-t)\over \ln(1+T)}\,.\eeq
In the following, a prime will denote differentiation w.r.t.~$t\in
[0,T]$, while the upper dot means a derivative w.r.t.~$s\in [0,1]$.
We claim that, by setting $\alpha\doteq\sqrt{\ln(1+T)}$ and defining
\bel{udef} u(t)\doteq u^*(\psi(t))+{\sqrt{2}\over\alpha^2}
\psi'(t)\sin(\alpha^3t)\cdot w^*(\psi(t))\,, \eeq the corresponding
solution $t\mapsto (q(t), p(t))$ of (\ref{DE1}), (\ref{indata})
satisfies the estimates (\ref{approx2}), provided that $|\bar p|$ is
small and  $T$  is sufficiently large. This will be proved in
several steps.

{\bf 1.}  It will be convenient to work with the variable
$s=\psi(t)\in [0,1]$, and derive an evolution equation for $q,p$ as
functions of $s$. By the definition of $\psi$ in (\ref{psidef}) it
follows \bel{psinv}t(s)~=~\psi^{-1}(s)~=~ (1+T)\Big( 1- e^{-s\cdot
\ln(1+T)}\Big), \eeq \bel{etadef} \qquad\quad \frac{ds}{dt}~=~
\psi'(t)~=~{1\over \ln(1+T)}\,\cdot{1\over 1+T-t}~=~ \frac{
e^{s\cdot \ln(1+T)}}{(1+T)\ln(1+T)}~\doteq~\eta(s)\,. \eeq In turn,
the functions
$$s\mapsto (\tilde q(s),\tilde p(s), \tilde u(s))
\doteq \Big(q(\psi^{-1}(s)),\,p(\psi^{-1}(s)),\,u(\psi^{-1}(s))\Big)
$$
satisfy
$$\frac{d}{ds}\left(\begin{array}{ccc} \tilde q\\\ \\\tilde p
\end{array}\right)
~=~\frac{1}{\eta(s)}\,\left[\left(\begin{array}{ccc} A \tilde p\\\
\\-\frac{1}{2}\tilde p^\dagger B
\tilde p \end{array}
\right)~+~\left(\begin{array}{ccc} K\\\ \\-\tilde p^\dagger C\end{array}
\right){du\over dt}~+~\left({d u\over dt}\right)^\dagger
\left(\begin{array}{ccc} 0 \\\ \\ D\end{array}\right)
{du\over dt}\right]\,.$$

Differentiating (\ref{udef}) and recalling that $\psi'=\eta$, we find
$$\begin{array}{rl}\displaystyle\frac{du}{dt}&=~\dot u^*(s)\,\eta(s)
+{\sqrt{2}\over\alpha^2} \psi''\,\sin(\alpha^3\psi^{-1}(s)) w^*(s)+
\sqrt{2}\,\alpha \,\eta(s)\, \cos(\alpha^3\psi^{-1}(s))\, w^*(s) \cr
&\cr &\qquad\qquad + {\sqrt{2}\over\alpha^2}
\sin(\alpha^3\psi^{-1}(s))\, \dot w^*(s)\,\eta^2(s).\end{array}
$$

Notice that (\ref{etadef}) yields
$$\psi''(t) ~=~{1\over \ln(1+T)} {1\over (1+T-t)^2} ~=~
\frac{ e^{2 s\cdot \ln(1+T)}}
{(1+T)^2\ln(1+T)}\,.$$
Putting together the above computations, we finally obtain
\bel{4.6}\begin{array}{l}\displaystyle
\frac{d}{ds}
\left(\begin{array}{ccc} \tilde q\\\ \\\tilde p \end{array}\right)~=~
\left(\begin{array}{ccc} A \tilde p\\\ \\-\frac{1}{2}\tilde p^\dagger B
\tilde p \end{array}
\right)\frac{1}{\eta(s)}~+~\left(\begin{array}{ccc} K\\\ \\-\tilde p^\dagger
 C\end{array}\right)(\dot u^*(s)+\phi_1(s)+\phi_2(s))\cr
\displaystyle ~+~\Big[\sqrt{2}\alpha\cos(\alpha^3\psi^{-1}(s))
\,w^*(s)+\zeta(s)\Big]^\dagger\left(\begin{array}{ccc} 0 \\\
\\ D\end{array}\right)\Big[\sqrt{2}\alpha\cos(\alpha^3\psi^{-1}(s))
\, w^*(s)+\zeta(s)\Big]\cdot{\eta(s)}\,,
\end{array}\eeq
where the functions  $\phi_1,\phi_2, \zeta$ are given respectively
by \bel{4.7}\left\{\begin{array}{l}\displaystyle \phi_1(s)~=
~\sqrt{2}\,\alpha\cos(\alpha^3\psi^{-1}(s))\, w^*(s)\,,\cr \cr
\displaystyle \phi_2(s)~= ~ {\sqrt{2}\over\alpha^2}
\sin(\alpha^3\psi^{-1}(s))\left(\frac{\dot
w^*(s)}{\ln(1+T)}+w^*(s)\right) \frac{e^{s\cdot\ln(1+T)}
}{1+T}\,,\cr \cr \displaystyle \zeta(s)~=~\dot
u^*(s)~+~{\sqrt{2}\over\alpha^2}
\sin(\alpha^3\psi^{-1}(s))\left(\frac{\dot w^*(s)}
{\ln(1+T)}+w^*(s)\right)\,\frac{e^{s\cdot\ln(1+T)}}{1+T}\,.
\end{array}\right.
\eeq
Before we derive the basic estimates, it is convenient to
introduce two more variables, namely
\bel{pstardef}
p^*(s)\doteq w^*(s)^\dagger D(q^*(s), u^*(s))
w^*(s)\,,\qquad\qquad {\bf p}(s)\doteq {\tilde p(s)\over \eta(s)}\,.
\eeq
We observe that
$${d\over ds} \bfp ~= ~{1\over\eta}\,{d\tilde p\over ds} -{\dot\eta
\over \eta^2} \tilde p ~=~{1\over\eta}\,{d\tilde p\over ds} -\alpha^2\bfp\,.$$
In term of ${\bf p}$, the system (\ref{4.6}) takes the form
\bel{finsyst}\begin{array}{l}\displaystyle
\frac{d}{ds}
\left(\begin{array}{ccc} \tilde q\\\ \\{\bf p} \end{array}\right)~=~
\left(\begin{array}{ccc} A {\bf p}\\\ \\-\frac{\eta(s)}{2}{\bf p}^\dagger B
{\bf p} \end{array}
\right)~+~\left(\begin{array}{ccc} K\\\ \\-\eta(s){\bf p}^\dagger
 C\end{array}\right)(\dot u^*(s)+\phi_1(s)+\phi_2(s))-\left(\begin{array}{ccc}
0 \\\
\\ \alpha^2{\bf p}\end{array}\right)\cr
 \cr
\displaystyle ~+~\Big[\sqrt{2}\alpha\cos(\alpha^3\psi^{-1}(s))
\,w^*(s)+\zeta(s)\Big]^\dagger\left(\begin{array}{ccc} 0 \\\
\\ D\end{array}\right)\Big[\sqrt{2}\alpha\cos(\alpha^3\psi^{-1}(s))
\, w^*(s)+\zeta(s)\Big]\,.
\end{array}\eeq
Notice that all functions $A,K,B,C,D$
here depend on $\tilde q, u^*$.

{\bf 2.} To help the reader, we give here a heuristic argument motivating
our key estimate.

By (\ref{etadef}) it follows
\bel{etasm}
0<\eta(s)\leq {1\over \ln(1+T)}\,.
\eeq
From the second equation in (\ref{finsyst}) one obtains
$$\begin{array}{rl}\displaystyle{d\over ds} {\bf p}(s) &\displaystyle
=~ - \alpha^2 {\bf p}(s) +2\alpha^2
\cos^2(\alpha^3\,\psi^{-1}(s))\,w^*(s)^\dagger D w^*(s)+ \O(1)\cdot
\alpha\cr &\cr &\approx ~-\alpha^2{\bf p}(s)
+\alpha^2\,w^*(s)^\dagger D w^*(s)+\O(1)\cdot \alpha \,.
\end{array}$$
Notice that last approximation follows from the fact that the
function $s\mapsto \cos^2(\alpha^3\psi^{-1}(s))$ is rapidly
oscillating and has average $1/2$.

Performing an integration by parts, the solution to the Cauchy problem
$$\dot P(s) ~=~-\alpha^2\cdot P(s)
+\alpha^2 \,w^*(s)^\dagger D(s) w^*(s)\,,\qquad\qquad P(0)=0,$$
can be written as
$$\begin{array}{rl}
P(s)&=~\displaystyle
\int_0^s e^{-\alpha^2(s-\sigma)} \alpha^2\Big(w^*(\sigma)^\dagger
D(\sigma) w^*(\sigma)\Big)\,d\sigma\cr
&\cr
&=~w^*(s)^\dagger
D(s) w^*(s)- e^{-\alpha^2 s}\Big(w^*(0)^\dagger
D(0) w^*(0)\Big)\cr
&\cr
&\qquad -\displaystyle \int_0^s e^{-\alpha^2(s-\sigma)} \cdot
\left[{d\over d\sigma}\Big(w^*(\sigma)^\dagger
D(\sigma) w^*(\sigma)\Big)\right]\,d\sigma\cr
&\cr
&=~w^*(s)^\dagger
D(s) w^*(s)+\O(1)\cdot \alpha^{-2}.\end{array}
$$
Since $\alpha=\sqrt{\ln(1+T)}\to\infty$ as $T\to\infty$, we thus
expect the convergence ${\bf p}(s)\to p^*(s)$ uniformly for $s\in
[0,1]$, where $p^*$ is the function introduced in (\ref{pstardef}).
  In turn, the first equation in (\ref{finsyst})
yields
$${d\over ds} \tilde q(s)~\approx~A p^*(s) + K \dot u^*(s)\,.$$
Indeed, in the computation of $K\phi$, the rapidly oscillating terms cancel out
in the limit.

As $T\to\infty$, we thus expect
$\tilde q(s)\to q^*(s)$ uniformly for $s\in [0,1]$.
Moreover,  by (\ref{etadef}) and (\ref{pstardef}),
$\tilde p(s) = \bfp(s)\eta(s)\to 0$ as $T\to\infty$.
The remaining steps of the proof will render  entirely rigorous
the above argument.
\v
{\bf 3.} In this section, for future use, we provide estimates
on two types of rapidly oscillating integrals.  In both cases the
key ingredient is an integration by parts. We assume that
the functions $h,\beta$ are $\C^2$ on the closed interval $[0,1]$,
with $h'(s)>0$.

First, multiplying and dividing by $h'(s)$ we compute
\bel{oscint1}
 \begin{array}{rl}
 &\displaystyle \int_0^\tau \cos(\alpha^3h(s)) \beta(s)\,ds\cr
 &\cr
&\quad =~\displaystyle \left(\int_0^\tau
\cos(\alpha^3h(r))\,h'(r)\,dr\right)\cdot {\beta(\tau)\over
h'(\tau)} -\int_0^\tau \left(\int_0^s \cos(\alpha^3
h(r))h'(r)\,dr\right) \cdot \left( {d\over ds} \,{\beta(s)\over
h'(s)}\right)\,ds \cr &\cr &\quad =~\displaystyle {
\sin(\alpha^3h(\tau))-\sin(\alpha^3h(0))\over \alpha^3}
\cdot{\beta(\tau)\over h'(\tau)}\cr &\cr &\displaystyle\qquad
-\int_0^\tau { \sin(\alpha^3h(s))-\sin(\alpha^3h(0))\over\alpha^3}
\cdot {\beta'(s) h'(s) - \beta(s) h''(s)\over [ h'(s)]^2}\,ds\,.
\end{array}
\eeq
Of course, an entirely similar estimate is valid replacing the cosine
with a sine function.
Next, by similar methods we compute
\bel{oscint2}
 \begin{array}{rl}
 &\displaystyle\int_0^\tau e^{-\alpha^2 (\tau-s)} 2\alpha^2 \cos^2(\alpha^3
h(s))
 \beta(s)\,ds\cr
 &\cr
 &
\displaystyle\qquad =~\int_0^\tau e^{-\alpha^2 (\tau-s)}\alpha^2
 \beta(s)\,ds+\int_0^\tau \alpha^2 e^{-\alpha^2 (\tau-s)} \Big(2
\cos^2(\alpha^3 h(s))-1\Big)
 h'(s) {\beta(s)
 \over h'(s)}\,ds \cr
 &\cr
&\qquad =~I_1+I_2\,.
\end{array}
\eeq
$$\begin{array}{rl}
 \displaystyle
 I_1&=~\displaystyle
 \alpha^2 e^{-\alpha^2\tau}\left( \int_0^\tau e^{\alpha^2 s}\, ds
\right) \beta(\tau) - \alpha^2 e^{-\alpha^2\tau} \int_0^\tau
\left(\int_0^s e^{\alpha^2r}
\,dr\right)\beta'(s)\,ds\cr
&\cr
&=~\displaystyle
(1-e^{-\alpha^2\tau}) \beta(\tau) - e^{-\alpha^2\tau} \int_0^\tau
\left( e^{\alpha^2s}-1\right) \beta'(s)\,ds\,,
\end{array}
$$
\bel{I1est}|I_1-\beta(\tau)|~\leq
~e^{-\alpha^2\tau}|\beta(0)|+{1\over\alpha^2}\,
\|\beta'\|_{\L^\infty}\,.\eeq

$$\begin{array}{rl}
 \displaystyle
 &I_2~=~\displaystyle
 \alpha^2 e^{-\alpha^2\tau}\Bigg\{\int_0^\tau \Big(2
\cos^2(\alpha^3h(s))-1\Big)
 h'(s) \,ds \cdot
 e^{\alpha^2 \tau}\,
{\beta(\tau)\over h'(\tau)} \cr
 &\cr
 &\quad\displaystyle - \int_0^\tau\left(
 \int_0^s \Big(2 \cos^2(\alpha^3 h(r))-1\Big)
 h'(r) \,dr \right)\cdot
 e^{\alpha^2 s}
\left( {\alpha^2\beta(s) +\beta'(s) \over h'(s)}
 -{\beta(s) h''(s)
 \over [h'(s)]^2}\right)\, ds\Bigg\}\cr
 &\cr
 &=~\displaystyle
 \alpha^2 e^{-\alpha^2\tau}\Bigg\{{\sin (2\alpha^3 h(\tau))-\sin(
2\alpha^3h(0))\over 2\alpha^3}
\cdot
 e^{\alpha^2 \tau}\,
{\beta(\tau)\over h'(\tau)} \cr
 &\cr
&\qquad \displaystyle - \int_0^\tau
 {\sin (2h(s))-\sin (2\alpha^3 h(0))\over 2\alpha^3}\cdot
 e^{\alpha^2 s}\cdot
\left( {\alpha^2\beta(s) +\beta'(s) \over h'(s)}
 -{\beta(s) h''(s)
 \over [h'(s)]^2}\right)
\, ds\Bigg\}\cr
\end{array}
$$
\bel{I2est}|I_2|~\leq~{1\over\alpha}\Big(
\alpha^2\|\beta\|_{\L^\infty} + \|\beta'\|_{\L^\infty}\Big)\cdot
\left\{ \min_{s\in [0,1]} h'(s)\right\}^{-1} + {1\over\alpha}\,
\|\beta\|_{\L^\infty}\left\|{h''\over [h']^2}\right\|_{\L^\infty}\,. \eeq

{\bf 4.} For a fixed $T>0$, consider the solution to the Cauchy
problem (\ref{finsyst}) with initial data \bel{ind0} \tilde
q(0)=\bar q,\qquad\qquad \bfp(0)=0\,.\eeq Its solution
$s\mapsto(\tilde q(s), \bfp(s))$ can be obtained as the fixed point
of a Picard transformation.  Namely, the transformation
$(q,p)\mapsto(\Q{(q,p)},\,\P{(q,p)})$ whose components are
\bel{Qdef}\begin{array}{rl} &\Q{(q,p)}(\tau)~=\displaystyle~\bar q
+\int_0^\tau A(q(s), u^*(s)) p(s)\,ds + \int_0^\tau K(q(s),
u^*(s))(\dot u^*(s)+\phi_2(s))\,ds\cr &\cr &\displaystyle \quad -
\int_0^\tau\Bigg[ K_q(q(s), u^*(s))\cdot \Big( A(q(s), u^*(s))p(s)
+K(q(s) , u^*(s)) ( \dot u^*(s)+\phi_1(s)+\phi_2(s))\Big)\cr &\cr
&\displaystyle \qquad + K_u(q(s), u^*(s)) \dot
u^*(s)\Bigg]\cdot\left( \int_0^s\phi_1(r)\,dr\right) \,ds ~+~
 K(q(\tau), u^*(\tau))\int_0^\tau \phi_1(s)\,ds\,, \end{array}
\eeq

\bel{Pdef}\begin{array}{rl}
&\!\P{(q,p)}(\tau)~=~\displaystyle~\int_0^\tau
e^{-\alpha^2(\tau-s)}\cdot 2\alpha^2 \cos^2(\alpha^3\psi^{-1}(s))\,
w^*(s)^\dagger D(q(s), u^*(s))w^*(s)\,ds \cr &\cr &
\displaystyle~-\int_0^\tau e^{-\alpha^2(\tau-s)}\cdot \eta(s)
p(s)^\dagger\left( B(q(s), u^*(s)){p(s)\over 2}+C(q(s), u^*(s))(\dot
u^*(s) +\phi_1(s)+\phi_2(s))\right)\,ds\cr &\cr
&~+\displaystyle\int_0^\tau e^{-\alpha^2(\tau-s)}\cdot\Bigg[
\sqrt{2}\alpha\cos(\alpha^3\psi^{-1}(s)) \,w^*(s)^\dagger D(q(s),
u^*(s)) \zeta(s)\cr &\cr & \qquad\qquad \displaystyle +
\zeta(s)^\dagger D(q(s),
u^*(s))\sqrt{2}\alpha\cos(\alpha^3\psi^{-1}(s)) \,
w^*(s)+\zeta(s)^\dagger D(q(s), u^*(s))\zeta(s)\Bigg]\,ds\,.
\end{array}
\eeq
Notice that the last two integral terms in (\ref{Qdef}) are obtained from
$$\int_0^\tau K(q(s), u^*(s))\,\phi_1(s)\,ds\,,$$
after an integration by parts.

On the family of couples of continuous functions $(q,p):[0,1]\mapsto
\R^{n+n}$ we consider the equivalent norm \bel{normdef1}
\Big\|(q,p)\Big\|_*\doteq \sup_{s\in[0,1]} ~ \max\left\{e^{-\rho s}
|q(s)|\,,~ {e^{-\rho s}\over \kappa}\,|p(s)| \right\}. \eeq We claim
that, if the constants $\rho,\kappa$ are chosen sufficiently large,
depending on the functions $A,B,C,D,K$ but not on $\alpha$, then the
Picard transformation $(\Q,\P)$ is a strict contraction w.r.t.~this
equivalent norm. Namely, \bel{pest1} \left\| \Big(\Q(q,p)-\Q(\hat
q,\hat p)\,,~\P(q,p)-\P(\hat q,\hat p)\Big)\right\|_* ~<~ {1\over
2}\, \Big\|(q-\hat q,\,p-\hat p)\Big\|_*\,.\eeq Moreover, we claim
that, as $\alpha=\sqrt{1+T}\to\infty$, one has \bel{pest2} \left\|
\Big(\Q(q^*,p^*), \P(q^*,p^*)\Big)-(q^*, p^*)\right\|_*~\to~ 0. \eeq
The two claims (\ref{pest1})-(\ref{pest2}) will be proved in the
next two sections.   In turn, they yield \bel{pest3} \left\| (\tilde
q,\bfp)-(q^*, p^*)\right\|_*~<~ 2\left\|\Big(\Q(q^*,p^*),
\P(q^*,p^*)\Big) -(q^*, p^*)\right\|_*~\to~ 0 \eeq as
$\alpha\to\infty$. {}From (\ref{pest3}), the conclusions in
(\ref{approx2}) will follow easily.

{\bf 5.}  In this step we establish  the strict contraction property
(\ref{pest1}).  As in Remark 3, it suffices to prove (\ref{pest1})
assuming that all functions $q,\hat q,p,\hat p$ take values within
some (possibly large) bounded set.

Assume that $\delta\doteq \Big\|(q-\hat q, \, p-\hat p)\Big\|_*$, so
that \bel{normd}|q(s)-\hat q(s)|\leq \delta e^{\rho
s}\,,\qquad\qquad |p(s)-\hat p(s)|\leq \delta \kappa \,e^{\rho
s}\qquad\qquad\forall s\in [0,1]. \eeq By (\ref{etadef}) and
(\ref{4.7}) we have \bel{pzest}0<\eta(s)\leq
{1\over\alpha^2}\,,\qquad |\phi_1(s)|\leq C_1\alpha\,,\qquad
|\phi_2(s)|\leq {C_1\over\alpha^2}\,,\qquad |\zeta(s)|\leq
C_1\,.\eeq Here and in the following, by $C_1,C_2,\ldots$ we denote
constants depending on the functions $A,B,C,D,K, u^*, w^*$, but not
on $\alpha\doteq \sqrt{\ln(1+T)}$. Applying (\ref{oscint1}) to the
case where  \bel{hbdef}h(s) = \psi^{-1}(s) = (1+T)
(1-e^{-s\cdot\ln(1+T)}), \qquad\qquad \beta(s)=\sqrt 2 \,\alpha\,
w^*(s),\eeq \bel{hder} h'(s)~=~ (1+T)\ln(1+T) e^{-s\cdot \ln(1+T)}
~\geq ~\ln(1+T)\,, \qquad\qquad {h''(s)\over [h'(s)]^2}
~\leq~1\,,\eeq one finds \bel{intphi1} \left|\int_0^\tau
\phi_1(s)\,ds\right| \leq {C_2\over \alpha^2} \eeq for all $\tau\in
[0,1]$.
 Recalling (\ref{Qdef}) and using (\ref{intphi1}) we obtain
\bel{Qest1} \Big|\Q(q,p)(s)-\Q(\hat q, \hat p)(s)\Big|~\leq
~\int_0^s C_3\Big( |p(s)-\hat p(s)| +|q(s)-\hat q(s)| \Big) \, ds\,.
\eeq By (\ref{etadef}) we have  $\eta(s)\leq \alpha^{-2}$.  From
(\ref{Pdef}) it thus follows \bel{Pest1} \Big|\P(q,p)(\tau)-\P(\hat
q, \hat p)(\tau)\Big|~\leq ~\int_0^\tau e^{-\alpha^2(\tau-s)}
C_4\Big( \alpha^2 |q(s)-\hat q(s)| +{1\over\alpha} |p(s)-\hat
p(s)|\Big)\, ds\,. \eeq

The bounds (\ref{normd}) and (\ref{Pest1}) imply \bel{Pest2}
\begin{array}{rl}\displaystyle {e^{-\rho \tau}\over \kappa}
\Big|\P(q,p)(\tau)-\P(\hat q, \hat p)(\tau)\Big| &\leq~\displaystyle
{e^{-\rho \tau}\over\kappa}\cdot \int_0^\tau e^{-\alpha^2(\tau-s)}
C_4\Big( \alpha^2 \delta e^{\rho s}  +{\kappa\over\alpha} \delta
e^{\rho s} \Big) \, ds\,,\cr &\cr
&<~\displaystyle{e^{-(\alpha^2+\rho)\tau}\over\kappa}\cdot
C_4\delta\,e^{(\alpha^2+\rho)\tau}\cdot{\alpha^2+{\kappa\over\alpha}\over
\alpha^2+\rho}~\leq~ {\delta\over 2}\,,
\end{array}
\eeq provided that $\kappa > 2C_4$ and $\alpha$ is suitably large.

In a similar way, the bounds (\ref{normd}) and (\ref{Qest1}) imply
\bel{Qest2}
\begin{array}{rl}e^{-\rho s}\Big|\Q(q,p)(\tau)-\Q(\hat q, \hat p)(\tau)\Big|
&\leq~ \displaystyle e^{-\rho s}\cdot \int_0^\tau C_3\Big(
\kappa\delta e^{\rho s} +\delta e^{\rho s} \Big) \, ds\,,\cr &\cr
&\leq\displaystyle~C_3\delta\cdot{\kappa+1\over\rho}~ \leq~
{\delta\over 2}\,,
\end{array}
\eeq provided that $\rho\geq 2C_3(\kappa+1)$.

 {\bf 6.} In this step we estimate the distance between
$(q^*,p^*)$ and the fixed point $(\tilde q, \bfp)$
of the transformation $(\Q,\P)$.
We recall  that $q^*$ satisfies
$$q^*(\tau) =\bar q
+\int_0^\tau A(q^*(s), u^*(s)) p^*(s)\,ds + \int_0^\tau K(q(s),
u^*(s))\,\dot u^*(s)\,ds\,,$$
with $p^*$ defined at (\ref{pstardef}).
Comparing this with  (\ref{Qdef}), we obtain
\bel{Qdiff1}\begin{array}{rl}
&\displaystyle
\Big|\Q(q^*,p^*)(\tau)-q^*(\tau)\Big|~\leq ~\left|\int_0^\tau K(q(s),
u^*(s))\phi_2(s)\,ds\right|\cr &\cr &\displaystyle \qquad +
\int_0^\tau\Bigg| K_q(q(s), u^*(s))\cdot \Big( A(q(s), u^*(s))p(s)
+K(q(s) , u^*(s)) ( \dot u^*(s)+\phi_1(s)+\phi_2(s))\Big)\cr &\cr
&\displaystyle \qquad + K_u(q(s), u^*(s)) \dot
u^*(s)\Bigg|\cdot\left| \int_0^s\phi_1(r)\,dr\right| \,ds ~+~ \Big|
K(q(\tau), u^*(\tau))\Big| \cdot \left|\int_0^\tau
\phi_1(s)\,ds\right|\,.
\end{array}
\eeq The definition of $\phi_1,\phi_2$ at (\ref{4.7}) implies
$$|\phi_1(s)|\leq C_5\,\alpha\,,\qquad\qquad |\phi_2(s)|\leq C_5\,
\alpha^{-2}.$$ Using the estimate (\ref{intphi1}) we thus obtain
\bel{Qdiff2}\Big|\Q(q^*,p^*)(\tau)-q^*(\tau)\Big|~\leq~C_6\,
\alpha^{-1}.\eeq

Next, comparing (\ref{pstardef}) with (\ref{Pdef}), we obtain
\bel{Pdiff1}\begin{array}{rl} \displaystyle
\left|p^*(\tau)-\P(q^*,p^*)(\tau)\right|&\leq~\displaystyle
\left|p^*(\tau)-\int_{0}^{\tau}e^{-\alpha^2(\tau-s)}2\alpha^2
\cos^2(\alpha^3\psi^{-1}(s)) \cdot p^*(s)ds \right|\cr &\cr
&\qquad\qquad ~\displaystyle + \int_{0}^\tau e^{-\alpha^2(\tau-s)}
C_7\,\alpha\,ds~\doteq ~ J_1+J_2\,.
\end{array}
\eeq A straightforward computation yields \bel{J2est} |J_2|~\leq~
C_7\,\alpha^{-1}.\eeq
To estimate $J_1$, we use (\ref{I1est})-(\ref{I2est}) with
$\beta(s) = p^*(s)$, $h(s) = \psi^{-1}(s)$.
By (\ref{psinv}), this implies
$$h'(s)~ = ~(1+T)\ln(1+T) e^{-s\cdot \ln(1+T)}~\geq ~\ln (1+T)
~=~\alpha^2,\qquad\qquad \left|h''(s)\over [h'(s)]^2\right|\leq 1\,.$$
Recalling that $p^*(s)=w^*(s)=0$ for $s\in [0,\ve_0]$, we
thus obtain
\bel{J1est} |J_1|~\leq~C_8\,\alpha^{-1}. \eeq

{\bf 7.}
By choosing $T$, and hence also $\alpha=\sqrt{\ln(1+T)}$, sufficiently large,
the difference between $(q^*, p^*)$ and the fixed point
$(\tilde q, \bfp)$
of the transformation $(\Q, \P)$ can thus be rendered arbitrarily small,
in the norm $\|\cdot \|_*$ introduced at (\ref{normdef1}).
Since the constant $\kappa$ is independent of $T$,
the norm $\|\cdot\|_*$ is uniformly equivalent to
the $\C^0$ norm.
This establishes the last two estimates in (\ref{approx2})
when $p(0)=\bar p = 0$.  By (\ref{etadef}) and (\ref{pstardef}),
we have
$$|\tilde p(s)|~=~|\bfp(s)|\,\eta(s)~\leq ~{|\bfp(s)|\over \ln(1+T)}\,.$$
Since $\bfp(s)\to p^*(s)$ as $T\to \infty$, uniformly
for $s\in [0,1]$, this implies the uniform convergence  $\tilde p(s)\to 0$. By continuity, all the estimates in (\ref{approx2})
remain valid whenever $|p(0)|\leq \delta$ for some $\delta>0$
small enough.
\endproof

\section{Proof of Theorem 2}
\label{sec:4}
\setcounter{equation}{0}
By a translation of coordinates, it is not restrictive to assume that
$(\bfq,\bfu)=0\in\R^{n+m}$.
By assumption,
when $\lambda=0\in\R^{n+m}$ we thus have
$(q^0(1),\,u^0(1))= (\bfq,\bfu) = 0\in\R^{n+m}$. Moreover
the $(n+m)\times (n+m)$ Jacobian matrix
\bel{Jmat}J~=~\left({\partial(q^\lambda(1),\,u^\lambda(1))\over
\partial\lambda}\right),\eeq
computed at the point $\lambda=0\in\R^{n+m}$, has maximum rank.
For notational convenience, we denote by $z=(q,u)$ the variable in $\R^{n+m}$
and
call $J^{-1}$ the inverse of the matrix $J$ in (\ref{Jmat}).
Taking $\lambda=J^{-1}z$, we thus have
\bel{5.2}\lim_{z\to 0} ~{z-\Big(q ^ {J^{-1}z}(1)~,~
u^{J^{-1}z}(1)\Big)\over |z|}~=~0.\eeq
Choosing $\rho>0$ sufficiently small, from
(\ref{5.2}) we deduce
\bel{5.3}\left|z-\Big(q ^ {J^{-1}z}(1)~,~
u^{J^{-1}z}(1)\Big)\right|~\leq~{|z|\over 3}\,,\qquad\hbox{for all}~~z\in
{\mathcal B}_\rho\,.\eeq
where ${\mathcal B}_\rho$ is the closed ball in $\R^{n+m}$,
centered at tho origin with radius $\rho$.

Next, we apply Lemma 1 and obtain a continuous map $(s,\lambda)\mapsto
(\tilde u^\lambda(s), \tilde w^\lambda(s))$ such that the corresponding
solutions
of (\ref{qapp}) satisfy (\ref{app5}) with $\ve = \rho/3$. Together with
(\ref{5.3}), this implies
\bel{5.4}\left|z-\Big(\tilde q ^ {J^{-1}z}(1)~,~
\tilde u^{J^{-1}z}(1)\Big)\right|~\leq~{2\rho\over 3}\,,\qquad\hbox{for
all}~~z\in
{\mathcal B}_\rho\,.\eeq

Finally, as in (\ref{udef}), we
define $\alpha=\sqrt{1+T}$ and the controls
\bel{Udef}
U^\lambda(t)\doteq \tilde u^\lambda(\psi(t))+\sqrt{2}\psi'(t)
\sin(\alpha\,t)\cdot \tilde w^\lambda(\psi(t))
\qquad\qquad t\in [0,T]\,.\eeq
If $|\bar p|$ is sufficiently small, choosing
$T$ sufficiently large the proof of Theorem 1
shows that the corresponding solutions
$t\mapsto\Big(Q^\lambda(t), P^\lambda(t)\Big)$
of (\ref{DE1})-(\ref{indata}) satisfy
\bel{5.6}
\left| \Big(Q^\lambda(\psi^{-1}(1)), \,U^\lambda(\psi^{-1}(1))\Big)
-\Big(\tilde q^\lambda(1),\,
\tilde u^\lambda(1)\Big)\right|\leq {\rho\over 3}\,,
\qquad\hbox{for all}~~\lambda\in
\Lambda\,.\eeq

We now consider the map
$$z~\stackrel{\Phi}{\longrightarrow}~
z~-~\Big( Q^{J^{-1}z}(1)~,~ U^ {J^{-1}z}(1)\Big)
\qquad\qquad z\in {\mathcal B}_\rho\,.$$
By (\ref{5.4}) and (\ref{5.6}), $\Phi$ is a continuous map of
the closed ball ${\mathcal B}_\rho$ into itself.  Hence, by Brouwer's
theorem, it has a fixed point $z^*$.
 This implies that exist   $\lambda^* = J^{-1}z^*$
such that  $ (Q^{\lambda^*}(1)~,~U^{\lambda^*}(1))~=~0\in\R^{n+m}$, completing
the proof. \endproof

\section{Examples}
\label{sec: A}
\setcounter{equation}{0}
In Lemma 1, the assumption (H) on the continuity of
the cone $\Gamma$ plays a key role. Indeed, if the map
$(q,u)\mapsto \Gamma(q,u)$ is not continuous the conclusion may be false.

{\bf Example 1.}  Let $q,u\in\R$, and consider the Cauchy
problem
\bel{ex1}\dot q= q^2\,\dot u^2\,,\qquad\qquad q(0)=-1.
\eeq
This corresponds to (\ref{DI1})-(\ref{Gammadef}), taking
$$K(q,u)\equiv 0,\qquad\quad A(q,u)= q^2\,,\quad\qquad D(q,u)\equiv 1.$$
In this case we have
$$\Gamma(q,u)= \left\{\begin{array}{rll}
&\{p\in\R\,;~ p\geq 0\}\qquad &\hbox{if}~~q\not= 0\,,\cr
&\{0\}\qquad &\hbox{if}~~q= 0\,,\cr
\end{array}\right.
$$
Hence, the map $\tilde q(t) = t-1$ provides a solution to
(\ref{DI1}).  However, for every $\C^1$ map $t\mapsto u(t)$
the corresponding solution of
(\ref{ex1}) satisfies
$q(t)<0$ for all $t\geq 0$.   Hence the map $\tilde q$ cannot be approximated
by smooth solutions of (\ref{ex1}).
\v
Next, we illustrate a simple application of Theorems 1 and 2.

{\bf Example 2.}
Consider a
 bead with mass $m$, sliding without friction along a bar.
 We assume that the bar can be rotated
around the origin on a horizontal plane (see fig.~1).   This system
can be described by two lagrangian parameters:
the distance $r$ of the bead from the origin, and the angle
$\theta$ formed by the bar and a fixed line through the origin.
The kinetic energy of the bead
is given by
\bel{bead1}T(r,\theta,\dot r,\dot\theta)= {m\over 2} (\dot r^2+r^2\dot
\theta^2)\,.
\eeq
We assign the angle $\theta=u(t)$ as
a function of time,
while the radius $r$ is the remaining free coordinate.
Setting $p=\partial T/\partial \dot r= m\dot r$,
the  motion is thus described  by
the equations
\bel{bead3}\left\{
\begin{array}{rl}
\dot r &=~p/m\,,\cr
\dot p &=~mr\dot u^2\,.\end{array}\right.
\eeq
Observe that in this case the right hand side of the equation
contains the square of the derivative of the control.

Consider the problem of steering the bead from
$A= (r_A, \theta_A)=(1,0)$ to a point
very close to $B=(r_B,\theta_B)= (1,\,\pi/2)$, during
an interval of time $[0,T]$ possibly very large.
Observe that this goal cannot be achieved by rotating the
bar with small but constant
angular velocity. Indeed,
choosing $\theta(t)=u(t)= \pi t/2T$, the trajectory of (\ref{bead3})
corresponding to the initial data $r(0)=1$, $p(0)=0$ is obtained by solving
$$\ddot r~=~\left(\pi\over 2T\right) ^2 r\,,\qquad\qquad r(0)=1,~~\dot
r(0)=0.$$
Hence
$r(t)={1\over 2} (e^{\pi t/2T} + e^{-\pi t/2T})$.
In particular,
$r(T)= {1\over 2} (e^{\pi/2}- e^{-\pi/2})$ for every choice of $T$.  Of course,
this value does not converge to
$1$ as $T\to\infty$.

We observe that,
in the present case, the differential inclusion  (\ref{DI1}) reduces to
$\dot r \geq 0$.
By Theorem 1, every continuous trajectory of the form
$t\mapsto (r(t), \theta(t))$, with $r$ a non-decreasing function of time,
can be tracked by solutions of the full system (\ref{bead3}).
In particular, according to (\ref{psidef}),
the trajectory
$$s\mapsto (r(s), \theta(s))= (1, \pi s/2)\qquad\qquad s\in [0,1]$$
can be traced by using the control
$$u(t) = \left( 1- {\ln(1+T-t)\over\ln(1+T)}\right) {\pi\over 2}\qquad\qquad
t\in [0,T].$$

Next, we observe that, if $r^*>r_0$, then
the point $(r^*,\theta^*)$
is normally reachable from the initial  point $( r_0, \theta_0)$ by solutions
of the differential
inclusion $(\dot r, \dot \theta)\in \R_+\times\R$.
Hence, by Theorem 2, for each $(r^*,\theta^*)$ with $r^*> r_0$ there exists
$T>0$ sufficiently large and
a control $u:[0,T]\mapsto\R$ with $u(0)=\theta_0$,
$u(T)=\theta^*$, such that
the solution of
(\ref{bead3}) with initial data
$$r(0)=r_0, \qquad\quad p(0)=0$$
satisfies $r(T)=r^*$.

\section{Derivation of the evolution equations}
\label{sec:B}
\setcounter{equation}{0}
Consider a system whose state is described by $N$ Lagrangian variables
$q^1,\ldots, q^N$.  Let the
 kinetic energy $T=T(q,\dot q)$
be given by a positive definite quadratic form
of the time derivatives $\dot q^i$, namely
\bel{3.1}
T(q,\dot q)~=~\frac{1}{2} \dot q^\dagger G \dot q~=~{1\over 2}
\sum_{i,j=1}^N g_{ij}(q) \,\dot q^i\dot q^j\,.
\eeq
Let the coordinates be split in two groups:
$\{q^1,\ldots, q^n\}$ and $\{ q^{n+1},\ldots, q^{n+m}\}$, with $N=n+m$.
The $(n+m)\times (n+m)$ symmetric
matrix $G$ in (\ref{3.1}) will thus take the corresponding block form
\bel{Gblock}
G=\left(\begin{array}{cc} G_{11}& G_{12}\\ G_{21} &G_{22}\end{array}
\right)=
\left(\begin{array}{cc} (g_{ij})& (g_{i, n+\beta})\\
(g_{n+\alpha, j})  &(g_{n+\alpha,\,n+\beta})
\end{array}
\right)\eeq
and
denote its inverse by
$$ \Hat G\doteq G^{-1}=\left(\begin{array}{cc} \Hat G_{11}& \Hat G_{12}\\
\Hat G_{21} &\Hat G_{22}\end{array}
\right)=
\left(\begin{array}{cc} (g^{ij})& (g^{i, n+\beta})\\
(g^{n+\alpha, j})  &(g^{n+\alpha,\,n+\beta})
\end{array}
\right).$$
Introduce the matrices
\bel{AKE}
A =\left(a^{ij}\right)\doteq (G_{11})^{-1}\,,\qquad\qquad
E=\left(e_{\alpha,\beta}\right)= (\Hat G_{22})^{-1},
\qquad\qquad
K= -A G_{12}\,,
\eeq
Since $\Hat G~=~G^{-1}$, we observe that
\bel{Ide}
\begin{array}{cc}\displaystyle \Hat{G}_{11}G_{11}+\Hat{G}_{12}G_{21}=Id\, , \qquad
\Hat{G}_{22}G_{22}+\Hat{G}_{21}G_{12}=Id\,,\cr\cr
\Hat{G}_{21}G_{11}+\Hat{G}_{22}G_{21}=0\, , \qquad G_{11}\Hat{G}_{12}+G_{12}\Hat{G}_{22}=0
\end{array}\eeq
Therefore, a straight forward rewriting of the above equations:
\bel{}
\left\{\begin{array}{cc}
G_{11}^{-1}=\Hat G_{11}+\Hat G_{12}G_{21}\Hat G_{11}^{-1} \, ,\qquad
\Hat G_{21}G_{11}^{-1} = - \Hat G_{22}^{-1}\Hat G_{21}\, ,
\cr\cr

\Hat G_{22}^{-1}= G_{22}+\Hat G_{22}^{-1}\Hat G_{21} G_{12} \, ,\qquad
\Hat G_{22}^{-1}\Hat G_{21} = -  G_{21} G_{11}^{-1}\,
\cr\cr

G_{11}\Hat G_{12}=-G_{12}\Hat G_{22} \, ,\qquad  A=G_{11}^{-1}\, ,\qquad   E=\Hat G_{22}^{-1}\end{array}\right.
\eeq

shows that the following identities hold
\bel{matrixid}
A=\Hat{G}_{11}- \Hat{G}_{12}E\Hat{G}_{21}\,,\qquad
E=G_{22}-G_{21}AG_{12}\,,\qquad
K=\Hat G_{12}E\,.\eeq
We assume that
a controller can prescribe the values
of the last $m$ coordinates
as functions of time, say
\bel{3.3}
q^{n+\alpha}(t)=u_\alpha(t)\qquad\qquad \alpha=1,\ldots,m\,,
\eeq
by implementing $m$ frictionless constraints.
Here {\bf frictionless} means that the forces produced by the constraints
make zero work in connection with any virtual displacement of the
remaining free coordinates $q^1,\ldots,q^n$.
In the absence of external forces, the motion is thus governed by the equations
\bel{lagreq}
\frac{d}{dt}
\frac{\partial T}{\partial \dot q^i} (q,\dot q)-\frac{\partial T}{\partial q^i}
(q,\dot q)~=~ \Phi_i(t)\qquad\qquad i=1,\ldots, n+m\,.
\eeq
Here $\Phi_i$ are the components of the forces generated by the constraints.
The assumption that these constraints are frictionless is expressed by the
identities
\bel{frictionless}
\Phi_i (t)\equiv 0\qquad\qquad i=1,\ldots,n\,.
\eeq
 Introducing the conjugate momenta
\bel{conjmom}
p_i~=~p_i(q,\dot q)~\doteq~ {\partial T\over\partial \dot q^i}
~=~ \sum_{i=1}^{n+m} g_{ij}(q)\,\dot q^j\,,
\eeq
We now consider
the system of Hamiltonian equations
for the first $n$ variables
\bel{redham}\left\{ \begin{array}{rl}\dot q^i&~=~~{\partial H\over\partial p_i}(q,p)\cr
\dot p_i&~=~-{\partial H\over\partial q^i}(q,p)\end{array}\right.
\qquad\qquad i=1,\ldots,n\,.
\eeq
Notice that (1.6) is a system of $2n$ equations for
$q^1,\ldots,q^n, p_1,\ldots,p_n$, where the right hand side
also depends on the remaining components $q^i,p_i$, $i=n+1,\ldots,n+m$.
We can remove this explicit dependence by inserting the values
\bel{closure}
\left\{\begin{array}{rl} q^{n+i}&=u_i(t)\,,\qquad\qquad  \dot q^{n+i}=\dot u_i(t)
\qquad\qquad \qquad  i=1,\ldots,m\,,\cr
p_j&=p_j(p_1,\ldots,p_n,\, \dot q^{n+1},\ldots, \dot q^{n+m})
\qquad\qquad j= n+1,\ldots,n+m\,.\end{array}\right.\eeq
 
From now on it will be more convenient to use vector notations. We thus write
$( q, u)=( q^1,\ldots,  q^n,  u_1,\ldots,
u_m)$,
$~ (p,\eta)= (p_1,\ldots, p_n, p_{n+1}, \ldots, p_{n+m})$.
Recalling that $\Hat G= G^{-1}$, we thus have
\bel{pqinv}
\left(\begin{array}{c}p\\   \eta\end{array}\right)
=\left(\begin{array}{cc} G_{11}& G_{12}\\ G_{21} &G_{22}\end{array}
\right)\left(\begin{array}{c}\dot q\\ \dot u\end{array}\right),
\qquad\qquad
\left(\begin{array}{c} \dot q\\ \dot u\end{array}\right)
=\left(\begin{array}{cc} \Hat G_{11}& \Hat G_{12}\\ \Hat G_{21} &\Hat G_{22}\end{array}
\right)\left(\begin{array}{c}p\\  \eta\end{array}\right).
\eeq

Multiplying by $A\doteq  (G_{11})^{-1} $ the identity
$$p=  G_{11} \dot q + G_{12} \dot u$$
we obtain
\bel{qdot1}\dot q~= ~(G_{11})^{-1} p - (G_{11})^{-1} G_{12} \dot u~=~
Ap- A G_{12} \,\dot u\,=Ap+K\,\dot u\,.
\eeq
Similarly, multiplying by $E=
 (\Hat G_{22})^{-1} $ the identity
$$\dot u=  \Hat G_{21} p  + \Hat G_{22} \eta $$
we obtain
\bel{eta} \eta~=~  (\Hat G_{22})^{-1} \dot u -(\Hat G_{22})^{-1}\Hat G_{21} p ~
=~  E \dot u - E \Hat  G_{21}\, p \,.
\eeq
{}From the equation
$$\dot p~=~ -\frac{1}{2} p^\dagger \frac{\partial \Hat G_{11}}{\partial q}
p
-p^\dagger \frac{\partial \Hat G_{12}}{\partial q} \eta
-\frac{1}{2} \eta^\dagger \frac{\partial \Hat G_{22}}{\partial q} \eta\,,$$
using (\ref{eta}) we obtain
\bel{pdot1}\dot p~=~ -\frac{1}{2} p^\dagger \frac{\partial \Hat G_{11}}{\partial q}p
-p^\dagger \frac{\partial \Hat G_{12}}{\partial q}  \Big(E \dot u - E \Hat  G_{21}\, p\Big)
-\frac{1}{2}
\Big(E \dot u - E \Hat  G_{21}\, p\Big)^\dagger \frac{\partial \Hat G_{22}}{\partial q}
\Big(E \dot u - E \Hat  G_{21}\, p\Big)\,.\eeq
Here and in the sequel, we use the notation $p^\dagger, \eta^\dagger\ldots $
to denote the transpose
of column vectors such as $p, \eta$.\vs
Now,we rewrite equation (6.15) into the form
\bel{dotp}
\begin{array}{lr}\displaystyle \dot{p}=-\frac{1}{2}p^{\dagger}\Big(\frac{\partial\Hat{G}_{11}}{\partial q}-\frac{\partial\Hat{G}_{12}}{\partial q}E\Hat{G}_{21}
+\frac{1}{2}(E\Hat{G}_{21})^{\dagger}\frac{\partial \Hat{G}_{22}}{\partial q}E\Hat{G}_{21}\Big)p\cr\cr
\displaystyle -p^{\dagger}\Big(+\frac{\partial\Hat{G}_{12}}{\partial q}E-(E\Hat{G}_{21})^{\dagger}\frac{\partial\Hat{G}_{22}}{\partial q}E\Big)\dot{u}-\frac{1}{2}\dot{u}^{\dagger}\Big(E^{\dagger}\frac{\partial\Hat{G}_{22}}{\partial q}E\Big)\dot{u}
\end{array}
\eeq
since $E=\Hat{G}_{22}^{-1}$,together with $\Hat{G}$ is symmetric,we have
$$ \frac{\partial\Hat{G}_{22}}{\partial q}E~=~-\Hat{G}_{22}\frac{\partial E}{\partial q}~~~~~~~;~~~~~~~\Hat{G}_{12}^{\dagger}=\Hat{G}_{21}$$
Hence equation (6.16) can be further rewrite as
\bel{dotp}
\begin{array}{lr}\displaystyle \dot{p}=-\frac{1}{2}p^{\dagger}\Big(\frac{\partial\Hat{G}_{11}}{\partial q}-\frac{\partial\Hat{G}_{12}}{\partial q}E\Hat{G}_{21}
-\frac{1}{2}\Hat{G}_{12}\frac{\partial E}{\partial q}\Hat{G}_{21}\Big)p\cr\cr
\displaystyle -p^{\dagger}\Big(\frac{\partial\Hat{G}_{12}}{\partial q}E
+(\Hat{G}_{12})\frac{\partial E}{\partial q}\Big)\dot{u}+\frac{1}{2}\dot{u}^{\dagger}\frac{\partial E}{\partial q}\dot{u}
\end{array}
\eeq
 Recall that $A=\Hat{G}_{11}-\Hat{G}_{12}E\Hat{G}_{21}$ and $K=\Hat{G}_{12}E$ , we have
$$\frac{\partial A}{\partial q}=\frac{\partial\Hat{G}_{11}}{\partial q}-
2\frac{\partial\Hat{G}_{12}}{\partial q}E\Hat{G}_{21}-\Hat{G}_{12}\frac{\partial E}{\partial q}\Hat{G}_{21}~~~~~~~~~~~~~~~~~~~~\frac{\partial K}{\partial q}=\frac{\partial\Hat{G}_{12}}{\partial q}E+\Hat{G}_{12}\frac{\partial E}{\partial q}.$$
Together with (6.13), we finally obtain that the evolution of the first $n$ variables
$(q^1,\ldots, q^n)$ and of the corresponding momenta
$(p_1,\ldots, p_n)$ can be described by the system
\bel{Smech*}
\left(\begin{array}{c}\dot q\\\,\\\dot p\end{array}\right) =
  \left(\begin{array}{c} Ap\\\, \\-\frac{1}{2}p^\dagger\frac{\partial
A}{\partial
q}p \end{array}\right)
+\left(\begin{array}{c} { K}\\\,\\
-p^\dagger\frac{\partial { K}}{\partial q} \end{array}\right)\dot u\, +\,\,
\dot u^{\dag} \left(\begin{array}{c} {0}\\\,\\
\frac{1}{2} \frac{\partial E}{\partial
q}\end{array}\right)\dot u  \,.
\eeq
Here $A,K,E$ are functions of $q,u$, defined as
\bel{AKE}
A =\left(a^{ij}\right)\doteq (G_{11})^{-1}\,,\qquad\qquad
E=G_{22}- G_{21}AG_{12}\,,
\qquad\qquad
K= -A G_{12}\,.
\eeq

\end{document}